\documentclass[12pt]{amsart}
\usepackage[left=3cm, right=3cm, top=3.5cm, bottom=4cm]{geometry} 
\usepackage[utf8]{inputenc}
\usepackage{amssymb}
\usepackage{amscd}
\usepackage{amsmath}
\usepackage{enumerate}
\usepackage{comment}
\usepackage{dsfont}
\usepackage{amsthm}
\usepackage{url}
\usepackage{ mathrsfs }
\usepackage{textcomp}
\usepackage[all]{xy}
\usepackage{tikz-cd}
\usepackage{xcolor}
\definecolor{darkgrey}{rgb}{0.4,0.4,0.5}
\usepackage{hyperref}
\hypersetup{
    colorlinks=true,
    linkcolor=blue,
    citecolor=darkgrey
}

\usepackage[OT2,T1]{fontenc}
  \DeclareSymbolFont{cyrletters}{OT2}{wncyr}{m}{n}
 \DeclareMathSymbol{\Sha}{\mathalpha}{cyrletters}{"58}

\usepackage{url}

\usepackage{multirow}

\newtheorem{theorem}{Theorem}[section]
\newtheorem{lemma}[theorem]{Lemma}

\newtheorem{conjecture}[theorem]{Conjecture}

\newtheorem{definition}[theorem]{Definition}
\newtheorem{proposition-definition}[theorem]{Proposition-Definition}

\newtheorem{question}[theorem]{Question}

\theoremstyle{remark}
\newtheorem{remark}[theorem]{Remark}

\numberwithin{equation}{section}

\newcommand\EatDot[1]{}

\newcommand{\cyc}{{\mathrm{cyc}}}

\newcommand{\Gal}{{\mathrm{Gal}}}
\newcommand{\Ker}{{\mathrm{Ker}}}

\newcommand{\Sel}{{\mathrm{Sel}}}

\newcommand{\Hom}{\mathrm{Hom}}

\newcommand{\ds}{\displaystyle}
\newcommand{\Ext}{\mathrm{Ext}}

\newcommand{\Sss}{S^\mathrm{ss}}
\newcommand{\Trace}{\mathrm{Trace}}

\usepackage{amsmath,mathpazo,amsthm, amssymb, mathtools}
\usepackage[english]{babel}
\usepackage{graphicx}
\usepackage{microtype}
\usepackage{cleveref}
\usepackage{enumerate}
\usepackage{stmaryrd}
\usepackage[none]{hyphenat}
\usepackage{mathpazo}
\newcommand{\ot}{\otimes}

\usepackage{mdframed}

\usepackage{cancel}
\newcommand{\T}{\text}
\newcommand{\vs}{\vspace{0.1 in}}

\newcommand{\cG}{\mathcal{G}} \newcommand{\cH}{\mathcal{H}}
 
 \newcommand{\cL}{\mathcal{L}}
 
\newcommand{\cO}{\mathcal{O}}

\usepackage{graphicx}

\newcommand{\bA}{\mathbb{A}} \newcommand{\bB}{\mathbb{B}}
\newcommand{\bC}{\mathbb{C}} 
\newcommand{\bD}{\mathbb{D}}
\newcommand{\bE}{\mathbb{E}} \newcommand{\bF}{\mathbb{F}}

 \newcommand{\bN}{\mathbb{N}}
 
\newcommand{\bQ}{\mathbb{Q}} 
\newcommand{\bR}{\mathbb{R}}

 \newcommand{\bZ}{\mathbb{Z}}

 \newcommand{\fp}{\mathfrak{p}}

\usepackage{multirow}
\usepackage{mdframed}
\usepackage{tikz-cd}

\begin{document}
\title[Iwasawa theory of supersingular elliptic curves]{Conjecture A and $\mu$-invariant for  Selmer groups of supersingular elliptic curves}

\author{Parham Hamidi}
\address{Department of Mathematics, The University of British Columbia}
\curraddr{Room 121, 1984 Mathematics Road\\
	Vancouver, BC\\
	Canada V6T 1Z2}
\email{phamidi@math.ubc.ca}

\author{Jishnu Ray}
\address{Department of Mathematics, The University of British Columbia}
\curraddr{Room 121, 1984 Mathematics Road\\
	Vancouver, BC\\
	Canada V6T 1Z2}
\email{jishnuray1992@gmail.com}

\thanks{This work is supported by grants from the University of British Columbia and PIMS-CNRS research grant}


\maketitle
\begin{abstract}

 Let $p$ be an odd prime and let $E$ be an elliptic curve defined over a number field $F$ with good reduction at primes above $p$. In this survey article, we give an overview of some of the important results proven for the fine Selmer group and the signed Selmer groups over cyclotomic towers as well as the signed Selmer groups over $\bZ_p^2$-extensions of an imaginary quadratic field where $p$ splits completely. We only discuss the algebraic aspects of these objects through Iwasawa theory. We also attempt to give some of the recent results implying the vanishing of the $\mu$-invariant under the hypothesis of Conjecture A. Moreover, we draw an analogy between the classical Selmer group in the ordinary reduction case and that of the signed Selmer groups of Kobayashi in the supersingular reduction case.  We highlight properties of signed Selmer groups (when $E$ has good supersingular reduction) which are completely analogous to the classical Selmer group (when $E$ has good ordinary reduction). In this survey paper, we do not present any proofs, however we have tried to give references of the discussed results for the interested reader.

\end{abstract}
\setcounter{tocdepth}{1}
\tableofcontents
\section*{Introduction}
Classical Iwasawa theory began from the work of Iwasawa who explored the growth of ideal class groups in towers of number fields. 
Iwasawa theory for elliptic curves, which shares many of its fundamental ideas with classical Iwasawa theory, deals with the arithmetic of elliptic curves over various infinite extensions. Among the main objects of study in Iwasawa theory are Selmer groups of the Galois representations attached to the elliptic curves and a prime $p$. Over the cyclotomic $\bZ_p$-extension of a number field, when the Galois representation (or equivalently, the elliptic curve) has good ordinary reduction over a prime $p$, then the Pontryagin dual of the Selmer group is conjectured to be a finitely generated torsion module over the  Iwasawa algebra of a $p$-adic Lie group of dimension $1$ (cf. Mazur's Conjecture \ref{conj:Mazure}). This Conjecture has been proved by Kato for elliptic curves over $\bQ$ and for the cyclotomic $\bZ_p$-extension $\bQ_{\cyc}$ of $\bQ$. 

\vs

Mazur's Conjecture can be generalized for arbitrary $p$-adic Lie extensions provided the elliptic curve has good ordinary reduction at $p$. Given a compact $p$-valued $p$-adic Lie group, the corresponding Iwasawa algebra admits a nice structure theory for finitely generated torsion modules (up to pseudo-isomorphism) (cf. \cite{CoatesSchneiderSujatha}). Therefore, the dual Selmer group enjoys important algebraic properties, and in the ordinary reduction case, it has been studied extensively.

\vs

When the elliptic curve has good supersingular reduction at primes over $p$, however, the Pontryagin dual of the Selmer group, even though still finitely generated, is not a torsion module over the corresponding Iwasawa algebra. In fact, it is shown to have positive rank \cite[Theorem 2.6]{CoatesSujatha_book} for the cyclotomic $\bZ_p$-extension and it is believed to be true more generally.

\vs

Nevertheless, this problem can be partially addressed if we focus our attention to certain subgroups of the Selmer groups. The signed Selmer groups, which were first introduced by Kobayashi (in \cite{Kobayashi}), are subgroups of Selmer groups and the emerging picture seems to be that the signed Selmer groups enjoy many of the properties in the supersingular case that are enjoyed by the Selmer groups in the ordinary case. This survey intends to give evidence for this claim. 

\vs

Over the cyclotomic $\bZ_p$-extension $\bQ_{\cyc}/\bQ$, the dual of the signed Selmer groups (often called plus/minus Selmer groups) are finitely generated torsion modules over the Iwasawa algebra of $\Gal(\bQ_{\cyc} /\bQ)$. 
In \cite{CoatesSujatha_fineSelmer}, Coates and Sujatha considered the notion of the \textit{fine} Selmer group, which is a subgroup of the signed Selmer groups. Therefore, the Pontryagin dual of fine Selmer is a finitely generated module over the the corresponding Iwasawa algebra. 
They conjectured that the $\mu$-invariant of the dual fine Selmer is zero; this is known as Conjecture A to Iwasawa theorists. 

\vs

Thus, we are left with four finitely generated torsion modules, viz. dual fine Selmer, dual plus and minus Selmer, dual of the torsion submodule of the Selmer. This article attempts to bring in light the connections between $\mu$-invariant of these torsion modules. 

\vs 

This article is organized as follows.
In Section \ref{Sec:def}, we define the Selmer group, the fine Selmer group and the signed Selmer groups over cyclotomic towers as well as the signed Selmer groups over a $\bZ_p^2$-extension of an imaginary quadratic field where $p$ splits completely.

\vs 

In Section \ref{Preliminaries}, we give some of the definitions that we will use in the rest of this article, such as the Iwasawa algebra, pseudo-nullity, and the Euler characteristic. We discuss the structure theorem and use it to define the $\mu$-invariant and the characteristic polynomial of a finitely generated module over a given Iwasawa algebra. Furthermore, we recall Conjecture A and discuss how it is related to Iwasawa $\mu$-invariant conjecture over cyclotomic extensions. We finish Section \ref{Preliminaries} by describing what is known as Iwasawa main conjecture which beautifully relates an algebraically defined object to an analytically defined one both coming from the arithmetic of number fields.

\vs

Then, in Section \ref{Sec:PropSelmer}, we recall some of the most important results proved so far for the signed Selmer groups over the cyclotomic $\bZ_p$-extension and $\bZ_p^2$-extensions. We note an exact sequence, due to Kobayashi \cite{Kobayashi}, connecting the fine Selmer group with that of the signed Selmer groups. This exact sequence can be thought of as a tool to transport information from the fine Selmer group (like Conjecture A or its dual  being torsion) to the signed Selmer groups. Moreover, we talk about how the Euler characteristic varies as we climb up the Iwasawa theoretic tower from the cyclotomic $\bZ_p$-extension to a $\bZ_p^2$-extension; this is a result of Lei and Sujatha \cite{LeiSujatha} (see Theorem \ref{Thm:LeiSujEuler}). 

\vs 

In Section \ref{Sec:ConjeAimplies}, over a cyclotomic $\bZ_p$-extension, when the elliptic curve has supersingular reduction at primes above $p$,  we explain how Conjecture A implies that the $\mu$-invariant attached to the torsion part of the dual Selmer group is zero. 
To show this, we use important results proved by Billot \cite{Billot} and then generalized by Wingberg \cite{Wingberg}.

\vs 

Our goal in Section \ref{Sec:pAdicHodge} is to give another view of Kobayashi's signed Selmer groups for the cyclotomic $\bZ_p$-extension. Lei, Loeffler and Zerbes \cite{LeiLoefflerZerbes},  \cite{LeiZerbes} have redefined Kobayashi's signed Selmer groups using Fontaine's $p$-adic Hodge theory and Fontaine' ring $\widetilde{\mathbb{E}}$ (cf. see \cite[Sec. 2.1]{LeiZerbes}). Recently, Scholze realized that the Fontaine's ring $\widetilde{\mathbb{E}}$ is the tilt of the perfectoid field $\widehat{\mathbb{Q}_p(\mu_{p^\infty})}$. In Section \ref{Sec:pAdicHodge}, 
we start by recalling the notion of perfectoid spaces of Scholze \cite{scholzeberkeley} which gives a geometric understanding of Fontaine's $p$-adic Hodge theory and Fontaine--Winterberger Theorem \cite{FontaineWintenberger}. Then, we use the language of perfectoids to give the definition of signed Selmer groups.  Since the titling construction of Scholze works for any general perfectoid field, this poses a natural open question of whether it is possible to define signed Selmer groups for more general $p$-adic Lie extensions using Scholze's construction of perfectoid fields.

\vs

Finally, in Section \ref{Sec:Future}, we turn our attention to elliptic curves with good ordinary reduction at $p$ and we list some properties of the dual Selmer group that are analogous to the case of signed Selmer groups in the supersingular reduction case. To finish, we pose two questions that remain open for elliptic curves with supersingular reduction concerning signed Selmer groups and their $\mu$-invariants, while their analogue holds for Selmer groups in the ordinary reduction case.

\section{Definitions of Selmer, signed Selmer, and fine Selmer groups}\label{Sec:def}
Let $p$ be an odd prime. Let  $F$ be a number field, $F^{\prime}$ be a subfield of $F$, and $E/F^{\prime}$ be an elliptic curve with  good reduction at all primes above $p$. Let $S$ be the (finite) set of primes of $F^\prime$ above $p$ and the primes where $E$ has bad reduction.
\subsection{Selmer group}
Let $F_{\cyc}$ be the cyclotomic $\bZ_p$-extension of $F$, and let $\Gamma$ denote $\Gal(F_{\cyc}/F)$ which is topologically isomorphic to $\bZ_p$. For each integer $n\geq 0$, let $F_n$ be the sub-extension of $F_{\cyc}$ such that $F_n$ is a cyclic extension of degree $p^n$ over $F$. The $p$-Selmer group over $F_n$ is defined by the sequence
\[
0 \rightarrow \Sel_p(E/F_n) \rightarrow  \ H^1(F_n, E_{p^\infty}) \rightarrow  \displaystyle\prod_{w}H^1(F_{n,w},E),
\]
where $w$ runs over all finite primes of $F_n$ and $F_{n,w}$ is the completion of $F_n$ at prime $w$. For $v \in S$, we let
\begin{align*}
J_v(E/F_n) :=  \displaystyle\bigoplus_{w \mid v} H^1(F_{n,w},E)(p)\cong
\displaystyle\bigoplus_{w \mid v} \frac{H^1 (F_{n,w},E_{p^\infty})}{E(F_{n,w})  \otimes \bQ_p/\bZ_p},
\end{align*}
where the isomorphism is due to the Kummer map (c.f. \cite[Section 1.6]{CoatesSujatha_book}). Then, we have an exact sequence\footnote{We refer the reader to \cite[Sections 1.7 and 2.2]{CoatesSujatha_book} for more details about the equivalent definitions of the Selmer group.}
\begin{equation}\label{def:Selmer with J_n}
0 \rightarrow \Sel_p(E/F_n) \rightarrow H^1(F_S/F_n, E_{p^\infty}) \xrightarrow{\lambda_n}  \displaystyle\bigoplus_{v \mid S}J_v(E/F_n).
\end{equation}
Here $F_S$ is the maximal extension of $F$ unramified outside $S$ and $\lambda_n$ is consists of localization maps. Here, the choice of $S$ implies that $F_\cyc \subset F_S$. 
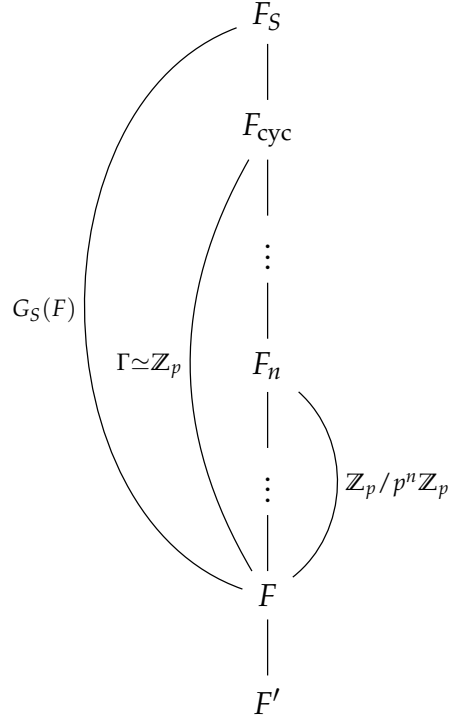
\[
\begin{tikzcd}
F_S \arrow[d, dash] \arrow[ddddd, dash, bend right=70, "G_S(F)" left]\\
F_\cyc \arrow[d, dash] \arrow[dddd,dash,  bend right=30, "\Gamma \simeq \bZ_p" left]\\
\vdots  \arrow[d, dash]\\
F_n \arrow[d, dash] \arrow[dd,dash,  bend left=50, "\bZ_p/p^n\bZ_p" right]\\
\vdots  \arrow[d, dash]\\
F \arrow[d, dash]\\
F^\prime 
\end{tikzcd}
\]
We define 
\[
\Sel_p(E/F_\cyc) := \varinjlim_n \Sel_p(E/F_n),
\]
and we obtain the following exact sequence by taking the direct limit of the exact sequence \eqref{def:Selmer with J_n} over intermediate field extensions $F_n$
\begin{equation*}\label{Selmercyc with J_v}
0 \rightarrow \Sel_p(E/F_\cyc) \rightarrow H^1(F_S/F_\cyc, E_{p^\infty}) \xrightarrow{\lambda_\cyc}  \displaystyle\bigoplus_{v \mid S}J_v(E/F_\cyc).
\end{equation*}

\begin{definition}
Suppose $K$ is a $p$-adic Lie extension of $F^\prime$. The $p$-Selmer group over $K$ is 
\[
\Sel_p(E/K) :=  \Ker \left(H^1(K, E_{p^\infty}) \rightarrow  \displaystyle\prod_{w}H^1(K_{w},E)\right),
\]
where $w$ runs over all finite primes of $K$ and $K_{w}$ is the union of the completions at $w$ of all finite extensions of $F^\prime$ contained in $K$.
\end{definition}
It is not hard to see that 
\begin{equation*}
    \Sel_p(E/K)=\varinjlim_{L}\Sel_p(E/L)
\end{equation*}
where $L$ runs over all finite extensions of $F^\prime$ contained in $K$.
\subsection{Signed Selmer groups over cyclotomic extension}\label{Sec:def:cyc}
We write $\Sss_p$ to be the set of primes of $F^\prime$ lying above $p$ where $E$ has supersingular reduction. Let $F^{\prime}_v$  be the completion of $F^{\prime}$ at a prime $v \in \Sss_p$ with the residue field $k^{\prime}_v$. Let $\widetilde{E}(k^{\prime}_v) $ be the $k^{\prime}_v$-points of the reduction of $E$ at place $v$. We assume the following:
\begin{enumerate}[(i)]
	\item $\Sss_p \not= \emptyset$;
\item For all $ v \in \Sss$ we have that $F^{\prime}_v$, the completion of $F^{\prime}$ at $v$, is $\bQ_p$;
	\item $a_v = 1+p -\#\widetilde{E}(k^{\prime}_v) =0$;
	\item $v$ is unramified in $F$.
\end{enumerate}

\vs

Let $\Sss_{p,F}$ denote the set of primes of $F$ above $\Sss$. Note that every ptime $v \in \Sss_{p,F}$ is totally ramified in the extension $F_\cyc/F$ and therefore there exists a unique prime in $F_n$ for each prime $v \in \Sss_{p,F}$. Let $F_{n,v}$ be the completion of $F_n$ at the unique prime over $v \in \Sss_{p,F}$. For every $v \in \Sss_{p,F}$, following Kobayashi \cite{Kobayashi} and Kitajima-Otsuki \cite{kitajima2018plus}, we define
\begin{align}
E^+(F_{n,v}) &= \{ P \in \hat{E}(F_{n,v}) \mid  \Trace_{n/m+1}P \in  \hat{E}(F_{m,v}) \text{ for all even } m, 0 \leq m \leq n-1 \},\label{eq:traceEven}\\
E^-(F_{n,v}) &= \{ P \in \hat{E}(F_{n,v}) \mid  \Trace_{n/m+1}P \in  \hat{E}(F_{m,v}) \text{ for all odd } m, -1 \leq m \leq n-1 \}\label{eq:traceOdd}.
\end{align}
Here $ \Trace_{n/m+1}$ is the trace map from   $\hat{E}(F_{n,v})$ to $\hat{E}(F_{m+1,v})$. 
 For all $n \geq 0$, we define 
\begin{align*}
E_{p^{n+1}}&=\Ker \left(  E(\overline{\bQ}) \xrightarrow{p^{n+1}} E(\overline{\bQ})  \right),\\
E_{p^{\infty}}&= \cup_{n \geq 0}E_{p^{n+1}}.    
\end{align*}

\vs

Identifying $E^{\pm}(F_{n,w})  \otimes \bQ_p/\bZ_p$ with a subgroup of $H^1(F_{n,w},E_{p^\infty})$ via the Kummer map, we may define the local terms 
\begin{align}
J_v^\pm(E/F_n) = \displaystyle\bigoplus_{w \mid v} \frac{H^1 (F_{n,w},E_{p^\infty})}{E^{\pm}(F_{n,w})  \otimes \bQ_p/\bZ_p}.
\end{align}
\begin{definition}
The plus and minus Selmer groups over $F_n$ are defined by 
$$\Sel_p^{\pm}(E/F_n) = \ker \left(   \Sel_p(E/F_n) \rightarrow \displaystyle\bigoplus_{v \in \Sss_{p,F}} J_v^\pm(E/F_n)  \right). $$
The plus and minus Selmer groups over $F_{\cyc}$ are obtained by taking direct limits
$$\Sel_p^\pm(E/F_\cyc) = \varinjlim_n \Sel_p^{\pm}(E/F_n).$$
\end{definition}
\subsection{Signed Selmer groups over a $\bZ_p^2$-extension}\label{sec:sec1.3}
Suppose $F$ is an imaginary quadratic field where $p$ splits completely. Let $F_\infty$ denote the compositum  of all $\bZ_p$-extensions of $F$. By Leopoldt's Conjecture we know that $G= \Gal(F_\infty/F)\simeq \bZ_p^2$, which implies that $F_\infty$ over $F_\cyc$ is a $\bZ_p$-extension. Let $v$ be a place of $F$ and $w$ be a place of $F_\infty$ above $v$. If $v \mid p$, then $F_v \simeq \bQ_p$ and $F_{\infty, w}$ is an abelian pro-$p$ extension over $F_v$ and $\Gal(F_{\infty,w}/F_v) \simeq \bZ_p^2.$ 
 Under this setting,  it is possible to define the plus and minus norm groups $E^{\pm}(F_{\infty,w}) \subset \hat{E}(F_{\infty,w})$ via Trace maps as in \eqref{eq:traceEven} and \eqref{eq:traceOdd}  (cf. Section 5.2 of \cite{LeiSujatha} which is a generalization of a construction by Kim \cite{Kim}). Similar to the above construction, we may define the local terms 
\begin{align}
J_v^\pm(E/F_\infty) = \displaystyle\bigoplus_{w \mid v} \frac{H^1 (F_{\infty,w},E_{p^\infty})}{E^{\pm}(F_{\infty,w})  \otimes \bQ_p/\bZ_p}.    
\end{align}
\begin{definition}
Retain the settings described above. The plus and minus Selmer groups over $\bZ_p^2$-extension $F_\infty$ are defined by 
 \begin{equation}
 \Sel_p^{\pm}(E/F_\infty) = \ker \left(   \Sel_p(E/F_\infty)  \rightarrow \displaystyle\bigoplus_{v \in \Sss_{p,F}}J_v^\pm(E/F_\infty)   \right).
 \end{equation}
\end{definition}
Here the (classical) Selmer group $\Sel_p(E/F_\infty)$ is defined by taking inductive limit of Selmer groups over all finite extensions contained in the $p$-adic Lie extension $F_\infty$.
\subsection{Fine Selmer groups}
The fine Selmer group of an elliptic curve over a number field is a subgroup of the Selmer group. The fine Selmer group is obtained when we impose extra stronger conditions at the primes above $p$ on the sequence defining the Selmer groups.
\begin{definition}\label{def:fineSelmer}
The \textbf{fine Selmer group} of $E/F^\prime$ over a number field $L/F^\prime$ is
\[
    \Sel_p^0(E/L):=\Ker \left(   H^1(L, E_{p^\infty}) \to \displaystyle\prod_\nu H^1(L_{\nu},E_{p^\infty})  \right),
\]
where $\nu$ runs over all places of $L$ and $L_{\nu}$ is the completion of $L$ at $\nu$. Let
\[
\Sel_p^0(E/F_\cyc):=\varinjlim_n \Sel_p^0(E/F_n).
\]
\end{definition}
When $K$ is a $p$-adic Lie extension of $F^\prime$, then 
\begin{equation*}
    \Sel^0_p(E/K)=\varinjlim_{L}\Sel^0_p(E/L)
\end{equation*}
where $L$ runs over all finite extensions of $F^\prime$ contained in $K$.
\begin{remark}
We can turn the direct product in the above definition \ref{def:fineSelmer} into a direct sum in the following way. Let $S$ be the set of primes of $F^\prime$ above $p$ and the primes where $E$ has bad reduction. Let $G_S(L)=\Gal(L_S/L)$ where $L_S$ denote the maximal extension of $L$ unramified outside of $S$ and the archimedean primes of $L$. Then, the above definition becomes (cf. \cite{CoatesSujatha_fineSelmer} for more details):
\[
    \Sel_p^0(E/L):=\Ker\left(  H^1(G_S(L), E_{p^\infty}) \to \displaystyle\bigoplus_{\nu \in S} K_\nu(L) \right),
\]
where
\[
    K_v(L)=\displaystyle\bigoplus_{\omega|\nu} H^1(L_\omega,E_{p^\infty}).
\]

\end{remark}

\section{Preliminaries of Iwasawa theory}\label{Preliminaries}
\subsection{Iwasawa Algebra}
\begin{definition}\label{Def:IwasawaAlgebra}\label{Sec:muInvariant}
Let $G$ as be a compact $p$-adic Lie group. We define the \textbf{Iwasawa algebra} of $ G $, denoted by $ \Lambda(G) $, to be the completed group algebra of $ G $ over $ \bZ_p $. That is  
\[
		\Lambda(G) =\bZ_{p}[[G]]:=\varprojlim_{N \subset G} \bZ_{p}[G/N],
\]
where $N$ runs over open normal subgroups of $G$.
\end{definition}
\begin{remark}
If there exists $n \geq 1 $ such that the $G\simeq \bZ_p^n$ (as topological groups), then $\Lambda(G)$ is isomorphic to the ring of formal power series $\bZ_p[[T_1,\cdots,T_n]]$, in indeterminate variables $T_1,\cdots,T_n$.
\end{remark}
When $G=\Gal(F_\cyc/F)\simeq\bZ_p$ (resp. $G\simeq\bZ_p^2$) then $\Lambda(G)\simeq\bZ_p[[T]]$ (resp. $\Lambda(G)\simeq\bZ_p[[T_1,T_2]]$). In this case, $\Lambda(G)$ is a commutative, regular local ring of Krull dimension $2$ (resp. $3$). More generally, when our Galois group $G$ is a non-commutative, torsion-free, compact $p$-adic Lie group of dimension $d$, then $\Lambda(G)$ is a non-commutative, left and right Noetherian local domain, and Auslander regular ring with global dimension $d+1$ (cf. \cite{venjakob2002structure}). This allows $\Lambda(G)$ to admit a nice dimension theory for its finitely generated modules\footnote{For information on the ring-theoretic properties of Iwasawa algebras we suggest Ardakov's survey paper \cite{ardakovring}.}.
\begin{definition}
\label{def:psuedonull}Let $M$ be a finitely generated $\Lambda(G)$-module. Then, $M$  is said to be \textbf{pseudo-null} if
\[
\dim_{\Lambda(G)}(M) \leq \dim_{\Lambda(G)}(\Lambda(G))-2.
\]
Furthermore, $M$ is $\Lambda(G)$-torsion if
\[
\dim_{\Lambda(G)}(M) \leq \dim_{\Lambda(G)}(\Lambda(G))-1.
\]
\end{definition}
Therefore, a pseudo-null $\Lambda(G)$-module is $\Lambda(G)$-torsion.
\begin{remark}
Equivalently, whenever $G$ is a torsion-free pro-$p$ group, a finitely generated $\Lambda(G)$-module $M$ is torsion if  Ext$_{\Lambda(G)}^0(M,\Lambda(G))=0$, and it is pseudo-null if \\Ext$_{\Lambda(G)}^i(M,\Lambda(G))=0$ for $i=0,1$.
\end{remark}
\begin{remark}
We note that if $\Lambda(G)$ is commutative, a $\Lambda(G)$-module $M$ is pseudo-null exactly when the annihilator ideal $\T{Ann}_{\Lambda(G)}(M)$  has height at least $2$. If $\Gamma:=\Gal(F_\cyc/F)\simeq\bZ_p$, and hence $\Lambda(\Gamma)\simeq\bZ_p[[T]]$, then we can show that a finitely generated $\Lambda(\Gamma)$-module is pseudo-null if and only if it is finite.
\end{remark}
\begin{remark}
For a general compact torsion-free $p$-adic Lie group $G$, a finite $\Lambda(G)$-module is always pseudo-null but the converse is not true.
\end{remark} 
\begin{definition}\label{pontryagin}
Suppose $G$ is a compact $p$-adic Lie group and $M$ is a $\Lambda(G)$-module. The \textbf{Pontryagin dual} of $M$, denoted by $\widehat{M}$, is defined to be 
\[
\widehat{M}:=\Hom(M,\bQ_p/\bZ_p),
\]
where $\Hom$ denotes the group of continuous homomorphisms.
\end{definition}
\subsection{Euler characteristic}
Suppose $G=\Gal(K/F)$, where $K$ is a $p$-adic Lie extension of a number field $F$. Suppose $G$ has finite dimension $d$ as a $p$-adic Lie group. Then, provided that $G$ has no non-trivial $p$-torsion, $G$ has $p$-cohomological dimension equal to $d$ (cf. Corollary 1 of \cite{serre1965dimension}). 

\begin{definition}\label{def:Eulerchar}
Suppose $G$ is a $p$-adic Lie group with $p$-cohomological dimension equal to $d$. Given a discrete primary $G$-module $M$ and if $H^i(G,M)$ is finite for $i=0,\cdots, d$, the \textbf{$G$-Euler characteristic}\footnote{For a more general see Definition 3.3.12 of \cite{Neukirch}.} of $M$ is defined to be
$$
\chi(G,M):=\ds\prod_{i=0}^{d}\left(\# H^i(G,M)\right)^{(-1)^i}.
$$
\end{definition}
\subsection{Structure Theorem}
 Suppose $K/F$ is a $\bZ_p$-extension, and therefore $G=\Gal(K/F)\simeq\bZ_p$. Moreover, suppose $M$ is a finitely generated $\Lambda(G)$-module. Then, the structure of $M$ as a $\Lambda(G)$-module is well understood up to pseudo-isomorphism. 
 
 \vs
 
 We say that a $\Lambda(G)$-homomorphism $f: M \rightarrow M^{\prime}$ between two $\Lambda(G)$-modules $M$ and $M^\prime$ is a pseudo-isomorphism if $f$ has finite kernel and cokernel. We say $M$ and $M^\prime$ are pseudo-isomorphic and we denote this by $M\sim M^\prime$. The following \textbf{structure Theorem} is due to Iwasawa.
\begin{theorem}\label{thm:structurethm}
Let $M$ be a finitely generated $\Lambda(G)$-module. Then 
\begin{align*}
    M \sim \Lambda(G)^r \oplus \left(  \displaystyle\bigoplus_{i=0}^s \Lambda(G)/p^{n_i}\Lambda(G) \right) \oplus \left(  \displaystyle\bigoplus_{j=1}^t  \Lambda(G)/f_j(T)^{m_j}\Lambda(G) \right),
\end{align*}
where $r$, $s$, $t$, $n_i$, and $m_j$ are non-negative integers and are unique up to reordering. Moreover, $f_j(T) \in \bZ_p[T]$ is an irreducible monic distinguished polynomial\footnote{A monic polynomial in $\bZ_p[T]$ is distinguished if $p$ divides all of its coefficients except for the leading one.} for each $j \in \{0,\cdots,t\}$.   
\end{theorem}
Note that in Theorem (\ref{thm:structurethm}), $r$ is the rank of $M$ as a $\Lambda(G)$-module and $M$ is torsion if and only if $r=0$. This enables us to define some important structural invariants for finitely generated $\Lambda(G)$-modules.

\begin{definition}
Given a $M$ be a finitely generated $\Lambda(G)$-module. Then the \textbf{characteristic polynomial} of $M$, denoted by $\T{charpoly}_{G}(M)$, is defined by
\[
\T{Charpoly}_{G}(M):=p^{\mu_G(M)}\ds\prod_{j=1}^t f_j(T)^{m_j},
\]
where $\mu_G(M):=\sum_{i=0}^s p^{n_i}$ and $s$, $t$, $n_i$, and $m_j$ are as in Theorem \ref{thm:structurethm}. Moreover, the \textbf{characteristic ideal} of $M$, denoted by $\T{Char}_{G}(M)$, is the principle $\Lambda(G)$-ideal generated by $\T{Charpoly}_{G}(M)$
\[
\T{Char}_{G}(M):=\left(\T{Charpoly}_G(M)\right)\Lambda(G)
\]
\end{definition}

It turns out that the above theorem can be naturally extended to any extension $K/F$ with $G=\Gal(K/F)\simeq \bZ_p^n$ where $n\geq 1$\footnote{ When  $G$ is a non-commutative torsion-free, compact $p$-adic Lie group, there is a structure theory for finitely generated torsion modules over the Iwasawa algebra $\Lambda(G)$  \cite{CoatesSchneiderSujathaII}.}. Suppose $G$ is a pro-$p$,
$p$-adic Lie group that is topologically isomorphic to $\bZ_p^n$ for some $n\geq 1$, and as a result, $\Lambda(G) \simeq \bZ_p[[T_1,\cdots,T_n]]$ for indeterminate variables $T_1,\cdots,T_n$. Suppose $M$ is a finitely generated $\Lambda(G)$-module. Recall from the Definition \ref{def:psuedonull} that we say $M$ is a pseudo-null if 
\[
\dim_{\Lambda(G)}(M) \leq \dim_{\Lambda(G)}(\Lambda(G))-2= (n+1)-2=n-1
\]
Moreover, we call a $\Lambda(G)$-homomorphism $f : M \to M^\prime$ between two finitely generated $\Lambda(G)$-modules a \textbf{pseudo-isomorphism} if the  kernel and cokernel of $f$ are pseudo-null $\Lambda(G)$-modules. 

\vs

Suppose $M$ is a finitely generated torsion $\Lambda(G)$-module. Let us denote the $p$-primary torsion $\Lambda(G)$-submodule of $M$ by $M(p)$. Then, we have the following Theorem \cite[Chapter VII, Section 4.4, Theorems 4 and 5]{Bourbaki4to7}:
\begin{theorem}\label{thm:formuinvariant}
Let $G$ be a pro-p
$p$-adic Lie group that is isomorphic to $\bZ_p^n$ for some $n\geq 1$ and let $M$ be a finitely generated torsion module over $\Lambda(G)$ with no elements of order $p$. Then, there exist a pseudo-isomorphism
\begin{equation*}
 M(p) \sim \ds\bigoplus_{i=0}^{s}\Lambda(G)/p^{n_i}\Lambda(G)   
\end{equation*}
where $n_1,\cdots, n_s$ are unique up to reordering.
\end{theorem}
The \textbf{$\mu$-invariant} of $M$, denoted by $\mu_G(M)$, is defined to be
\[
\mu_G(M)=\sum_{i=0}^s n_i,
\]
where $n_1,\cdots,n_s$ are as in
Theorem \ref{thm:formuinvariant} described above. 
\subsection{Introduction to Conjecture A}
Let $F$ be a number field and let $p$ be an odd prime. Let $F_\cyc/F$ be the cyclotomic $\bZ_p$-extension of $F$ and let $\Gamma:=\Gal(F_\cyc/F)$. Then, let $F_n/F$ be the sub-extension of $F_\cyc$ such that
\[
\Gal(F_n/F)=\Gamma/\Gamma^{p^n}=\bZ/p^n\bZ.
\]
Let $L_n/F_n$ be the maximal unramified abelian $p$-extension of $F_n$ (so $L_n$ is the $p$-Hilbert class field extension of $F_n$).
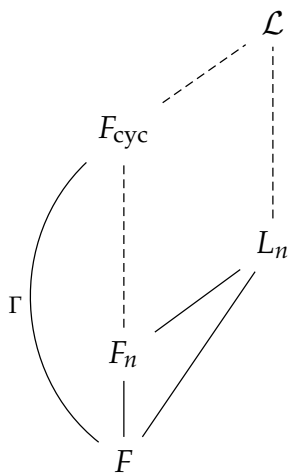
\[
\begin{tikzcd}
                                                                          & \cL                         \\
F_\cyc \arrow[ru, no head, dashed] \arrow[ddd, "\Gamma"', no head, bend right=49] &                                 \\
                                                                          & L_n \arrow[uu, no head, dashed] \\
F_n \arrow[ru, no head] \arrow[uu, no head, dashed]                       &                                 \\
F \arrow[u, no head] \arrow[ruu, no head]                                 &                                
\end{tikzcd}
\]
Note that by class field theory $\Gal(L_n/F_n)$ is isomorphic to the $p$-Sylow subgroup of the ideal class group of $F_n$. Let $\cL:=\cup_n L_n$, and then
\[
X(F_\cyc):=\Gal(\cL/F_\cyc)
\]
can be given a (natural) $\Gamma$-action, and therefore it can be turned into a $\Lambda(\Gamma)$-module. Iwasawa proved that $X(F_\cyc)$ is a finitely generated torsion $\Lambda(\Gamma)$-module (cf. \cite{IW-Gamma-extensions}). Moreover, in his celebrated work \cite{Iwasawaconjecture}, Iwasawa conjectured the following.
\vspace{0.1 in}

\textbf{Iwasawa $\mu$-invariant conjecture for cyclotomic extensions:} For any number field $F$, $X(F_\cyc)$ is a finitely generated $\bZ_p$-module. 

\vspace{0.1 in}
By structure theorem, this equivalent to say that $\mu_\Gamma(X(F_\cyc))=0$. We have the following important theorem of Ferrero--Washington \cite{Ferrero-Washington}.
\begin{theorem}\label{thm:Ferrero-Washington}
Iwasawa $\mu$-invariant conjecture holds for all abelian number fields.
\end{theorem}

\vs

Let $\mathfrak{X}^0(E/F_\cyc)$ denote the Pontryagin dual of $\Sel_p^0(E/F_\cyc)$. Then, we have the following conjecture of Coates and Sujatha \cite{CoatesSujatha_fineSelmer}.

\vspace{0.1 in}

\textbf{Conjecture A:} For any number field $F$, $\mathfrak{X}^0(E/F_\cyc)$ is a finitely generated $\bZ_p$-module. 

\vspace{0.1 in}

We note that this is equivalent to $\mathfrak{X}^0(E/F_\cyc)$ being a torsion $\Lambda(\Gamma)$-module and having $\mu$-invariant equal to zero (for the definition of $\mu$-invariant see Section \ref{Sec:muInvariant}).
\begin{remark}
See Section \ref{Numerical examples} for several examples of elliptic curves where it is shown that Conjecture A is satisfied. 
\end{remark}
The following theorem of Coates--Sujatha \cite[Theorem 3.4]{CoatesSujatha_fineSelmer} relates Conjecture A and the Iwasawa $\mu$-invariant conjecture for cyclotomic field extensions of algebraic number fields. 
\begin{theorem}\label{thm:equivalenceofconjAandIwaswawa}
Suppose $E/F$ is an elliptic curve over a number field $F$. Suppose $p$ is an odd prime such that the field extension $F(E_{p^\infty})/F$ is a pro-$p$ extension. Then Conjecture A holds for $E$ over $F_\cyc$ if and only if the Iwasawa $\mu$-invariant conjecture holds for $F_\cyc$.
\end{theorem}
Using Theorem \ref{thm:equivalenceofconjAandIwaswawa}, Coates and Sujatha gave the following conditions for an elliptic curve $E/F$ to satisfy Conjecture A \cite[Corollary 3.5]{CoatesSujatha_fineSelmer}.
\begin{theorem}
Suppose $p$ is an odd prime. If there exists a finite extension $L/F$ of $F$ such that
\begin{enumerate}

    \item $L \subset F(E_{p^\infty})$;
    \item $F(E_{p^\infty})$ is a pro-$p$, and
    \item the Iwasawa $\mu$-invariant conjecture hold for $L_\cyc$.
\end{enumerate}
Then Conjecture A holds for $E$ over $F_\cyc$. That is, $\mathfrak{X}^0(E/F_\cyc)$ is a finitely generated $\bZ_p$-module.
\end{theorem}
Moreover, Coates and Sujatha combined the above Theorem and Theorem \ref{thm:Ferrero-Washington} to prove \cite[Corollary 3.6]{CoatesSujatha_fineSelmer}:
\begin{theorem}
Suppose $p$ is an odd prime and let $F/\bQ$ be a number field such that $\Gal(F/\bQ)$ is abelian. If $E_{p^\infty}(F)\neq 0$ then Conjecture A holds for $E$ over $F_\cyc$. 
\end{theorem}
\subsection{Iwasawa main conjecture}
In this Section, we give a brief account of the classical Iwasawa Main Conjecture (or IMC for short). IMC describes a striking relation between two objects attached to number fields, a $p$-adic analytic object, in the form of the values of $p$-adic $L$-functions and certain algebraically defined Iwasawa module. Therefore, Iwasawa main conjecture gives us a bridge between algebraic treatments of arithmetic objects and that of the analytic side. To give the statement of Iwasawa main conjecture, we need to introduce more notations. We follow Chapter 15 of Washington's book \cite{washington1997introduction}.

\vs

Let $p$ be an odd prime. Let $F=\bQ(\mu_{p})$ and consider the cyclotomic extension $\bQ(\mu_{p^\infty})/\bQ(\mu_{p})=F_\cyc/F$, where $\mu_{p^n}$ denotes the $p^n$-th primitive root of unity and  $\bQ(\mu_{p^\infty})=\cup_{n\geq0}\bQ(\mu_{p^n})$. Let $\Gamma=\Gal(\bQ(\mu_{p^\infty})/\bQ(\mu_{p}))$. Furthermore, 
\begin{enumerate}
    \item Let $A_n$ be the $p$-Sylow part of the ideal class group of $\bQ(\mu_{p^{n+1}})$ for $n\geq0$ and $A_\infty:=\varinjlim_n A_n$.
    \item Let $\cL$ be the maximal unramified abelian $p$-extension of $\bQ(\mu_{p^\infty})$ and $X(F_\cyc)=\Gal(\cL/\bQ(\mu_{p^\infty}))$.
    \item Let $\varepsilon_i$ denote the $i^{th}$ idempotent element of $\Gal(\bQ(\mu_p)/\bQ)$, where $i$ is odd.
    \item Let $L_p(s,\omega^j)$ be the $p$-adic $L$ function for $\omega^j$ where $\omega$ is the Teichmüller character at $p$ and $j\neq$ is even (cf. Chapter 5-7 of \cite{washington1997introduction} for a construction $p$-adic $L$ function).
    \item Let $f(T,\omega^j) \in \bZ_p[[T]]\cong\bZ_p[[\Gamma]]$ be such that
    $f(1+p)^s-1,\omega^j)=L_p(s,\omega^j)$ (cf. Chapter 13, Section 6 of \cite{washington1997introduction}).
\end{enumerate}
Then, as we discussed, $X(F_\cyc)$ is a $\Lambda(\Gamma)$-module and Iwasawa proved that $X(F_\cyc)$ is finitely generated and torsion. Since $F=\bQ(\mu_p)$ is an abelian extension of $\bQ$. Theorem \ref{thm:Ferrero-Washington} implies that $\mu_{\Gamma}(X(F_\cyc))=0$ and by the structure Theorem \ref{thm:structurethm} we see that
\[
X(F_\cyc) \sim \ds \bigoplus_{j=1}^t \Lambda(\Gamma)\big/f_j(T)^{m_j} \Lambda(\Gamma)
\]
and 
\[
\T{Charpoly}_{\Gamma}(X(F_\cyc))=\ds\prod_{j=1}^tf_j(T)^{m_j} 
\]
where $f_j \in \Lambda(\Gamma)\cong \bZ_p[[T]]$ is an irreducible distinguished polynomial for each $j\in\{ 1,\cdots,t\}$, as described in structure Theorem (\ref{thm:structurethm}). The (classical) Iwasawa main conjecture was proved by Mazur and Wiles in \cite{Mazur-Wiles} and it 
states that (cf. \cite[Chapter 15, Section 4, Theorem 15.14]{washington1997introduction}):
\begin{theorem}
Retain the notations (1) to (5) described above. Let $i$ be an odd integer such that $i \neq 1 \pmod{p-1}$. Then
\[
\T{Charpoly}_{\Gamma}(\varepsilon_iX(F_\cyc))=f(T,\omega^{1-i}) u(T)
\]
where $u(T) \in \Lambda(\Gamma)^\times$ is a unit of $\Lambda(\Gamma)$. 
\end{theorem}
The Iwasawa main conjecture was generalized (and proved in many cases) to more general number fields and elliptic curves.

\section{Properties of fine Selmer and signed Selmer groups}\label{Sec:PropSelmer}
Let $\mathfrak{X}(E/F_\square)$ denote the Pontryagin dual of $\Sel_p(E/F_\square)$, where $\square\in\{\cyc,\infty\}$. Then $\mathfrak{X}(E/F_\square)$ is finitely generated as a $\Lambda(G)$-module, where $G=\Gal(F_\cyc/F)$ if $\square = \cyc$, and $G=\Gal(F_\infty/F)$ if $\square = \infty$. But unlike the ordinary case\footnote{See Conjecture \ref{conj:Mazure} for more details.}, it is believed (shown for the cyclotomic $\bZ_p$-extensions, see \cite[Theorem 2.5]{CoatesSujatha_book}) to have positive rank, and therefore it is no longer a torsion Iwasawa module. However, the Pontryagin dual of the signed Selmer groups $\Sel_p^{\pm}(E/F_\square)$, which we denote by $\mathfrak{X}^{\pm}(E/F_\square)$, is conjectured to be torsion.
\subsection{Cyclotomic case}
Let $E$ be an elliptic curve over $\bQ$ with good supersingular reduction at $p$ and suppose $a_p=0$. Let us denote $\Gal(\bQ_\cyc/\bQ)$ by $\Gamma$. Then, Kobayashi proved that (\cite{Kobayashi}, Theorem 1.2):
\begin{theorem}The  dual Selmer group $\mathfrak{X}^{\pm}(E/\bQ_\cyc)$ is
 a finitely generated torsion $\bZ_p[[\Gamma]]$-module.
\end{theorem}
Moreover, Kobayashi used Kato's Euler system to prove the following theorem \cite[Theorem 1.3]{Kobayashi}:
\begin{theorem}
Suppose $E/\bQ$ does not have complex multiplication. Then for almost all primes $p$, we have that the ideal generated by Pollack's $p$-adic $L$-function $\cL_p^\pm(E,X)$ as $\Lambda(\Gamma)$-modules (see \cite{pollack2003p} for more details) is contained in the characteristic ideal of the Pontryagin dual of $\Sel_p^{\pm}(E/\bQ_\cyc)$. That is,
\[
  (\cL_p^\pm(E,X)) \subseteq \T{Char}_\Gamma(\mathfrak{X}^{\pm}(E/\bQ_\cyc)).
\]
\end{theorem}
The Iwasawa main conjecture predicts that this containment is an equality. 

\vs

More generally, let us adapt the notations of Section \ref{Sec:def:cyc}. Here, $F^\prime$ is a number field and $E/F^\prime$ is an elliptic curve with good reduction at all primes above $p$, where $p$ is an odd prime. Suppose $F/F^\prime$ is a finite extension of $F^\prime$. Then, we we have the exact sequence
\begin{equation}\label{eq:lambdaplusminus} 
0 \to \Sel_p^\pm(E/F_\cyc) \to  H^1(F_S/F_\cyc,E_{p^\infty}) \xrightarrow{\lambda_\cyc^\pm} \displaystyle\oplus_{\nu\in S} J_v^\pm(E/F_\cyc). 
\end{equation}
Lei and Sujatha proved the following statement (\cite{LeiSujatha}, Proposition 4.4).
\begin{theorem}\label{thm:surjective}
$\mathfrak{X}^{\pm}(E/F_\cyc)$ is torsion as a $\Lambda(G)$-module if and only if $\lambda_\cyc^\pm$ is surjective in \eqref{eq:lambdaplusminus} and $H^2(F_S/F_\cyc, E_{p^\infty})=0$.
\end{theorem}
\begin{remark}
$H^2(F_S/F_\cyc, E_{p^\infty})=0$ is called the \textbf{weak Leopoldt's Conjecture} for $E_{p^\infty}$ over $F_\cyc$. This is equivalent to say that $\mathfrak{X}^0(E/F_\cyc)$, the Pontryagin dual of $\Sel_p^0(E/F_\cyc)$, is a finitely generated torsion $\Lambda(\Gamma)$-module \cite[Lemma 3.1]{CoatesSujatha_fineSelmer}. More generally, the weak Leopoldt's Conjecture for $E_{p^\infty}$ on a number field $F$ predicts that 
\[
H^2(F_S/F, E_{p^\infty})=0.
\]
\end{remark}
\begin{remark}
Note that if $\Sel_p(E/F)$ is finite then $\Sel_p^\pm(E/F)$ is finite and the dual signed Selmer groups $\mathfrak{X}^\pm(E/F_\cyc)$ are $\Lambda(\Gamma)$-torsion (cf. Remark 4.5 of \cite{LeiSujatha}). This provides examples for $\mathfrak{X}^{\pm}(E/F_\cyc)$ to be torsion. By the above Theorem, it also provides examples for the weak Leopoldt's Conjecture for $E_{p^\infty}$ over $F_\cyc$.
\end{remark}
For an elliptic curve $E$ defined over $\bQ$, we have:
\begin{theorem}[\cite{Kobayashi}, Theorem 5.1 and Corollary 7.2]\label{Kobayashi 7.2}
$\mathfrak{X}^0(E/\bQ_\cyc)$ is a finitely generated torsion $\Lambda(\Gamma)$-module.
\end{theorem}
In the following, we give an exact sequence due to Kobayashi \cite{Kobayashi} which connects the signed Selmer groups with that of the fine Selmer group.\\

Let $F_n=\bQ(\zeta_{p^{n+1}})$. Let $k_n$ denote $\bQ_p(\zeta_{p^{n+1}})$. Let $T_p(E)$ be the $p$-adic Tate module of elliptic curve $E$ and $V_p(E):=T_p(E) \otimes \bQ_p$. Suppose $S_n$ is the set of places of $F_n$ over $S$. $G_{n,S}:=\Gal(F_{n,S}/F_n)$ where $F_{n,S}$ is the maximal unramified extension of $F_n$ outside $S_n$. We let
\begin{align*}
    H^i_{/S}(T_p(E)):=\varprojlim_n H^i(G_{n,S},T_p(E)), \ \ H^i_{\T{Iw},\nu}(T_p(E)):=\varprojlim_n H^i(F_{n,\nu},T_p(E)).
\end{align*}
Moreover, we define $H_{\T{Iw},\pm}^1(T_p(E))$ as the exact annihilator with respect to the Tate pairing 
\begin{align*}
    H^1(k_n, V_p(E)/T_p(E)) \times H^1(k_n,T_p(E)) \to \bQ_p/\bZ_p
\end{align*}
of the subgroup $E^\pm(k_n) \otimes \bQ_p/\bZ_p \subseteq H^1(k_n,V_p(E)/T_p(E))$ (\cite{Kobayashi}). Additionally, we have the following exact sequence relating fine Selmer with the signed Selmer groups (\cite{Kobayashi}, Theorem 7.3):
\begin{align}
    0 \to H^1_{/S}(T_p(E)) \xrightarrow{\alpha} H_{\T{Iw}}^1(T_p(E))/H_{\T{Iw},\pm}^1(T_p(E)) \to \mathfrak{X}^\pm(E/\bQ_\cyc) \to \mathfrak{X}^0(E/\bQ_\cyc) \to 0.\label{eq:relating fine Selmer}
\end{align}
We note that if $T_p(E)/pT_p(E)$ is irreducible as a two dimensional representation of $Gal(\overline{\bQ}/\bQ)$ over $\bF_p$, then $H^1_{/S}(T_p(E))$ is a free $\Lambda(\Gamma)$-module of rank $1$ (\cite{Kato}, Theorem 12.4).
\begin{remark}
Kobayashi uses the Coleman maps and Kato's Euler system to show that the map $\alpha$ in (\ref{eq:relating fine Selmer}) is a non-zero map and the cokernel of $\alpha$ is a $\Lambda(\Gamma)$-torsion module (Theorem 6.3 of \cite{Kobayashi}). Then Theorem \ref{Kobayashi 7.2} and the exact sequence (\ref{eq:relating fine Selmer}) give a proof of the fact that $\mathfrak{X}^\pm(E/\bQ_\cyc)$ is torsion.
\end{remark}
The following Theorem tells us about the triviality of finite $\Lambda(\Gamma)$-submodules of dual signed Selmer groups of an elliptic curve $E/\bQ$ over the cylcotomic $\bZ_p$-extension (\cite{kimdu2013plus}, Theorem 1.1).
\begin{theorem}\label{sspsuedonull}
Let $p$ be an odd prime and suppose $F$ is a finite extension of $\bQ$ in which $p$ is unramified. Let $E/\bQ$ have good supersingular reduction at $p$  with $a_p=0$.
\begin{enumerate}
    \item Suppose $\mathfrak{X}^-(E/F_\cyc)$ is $\Lambda(\Gamma)$-torsion. Then, $\mathfrak{X}^-(E/F_\cyc)$ has no nontrivial finite $\Lambda(\Gamma)$-submodules.
    \item Suppose $\mathfrak{X}^+(E/F_\cyc)$ is $\Lambda(\Gamma)$-torsion, and that $p$ splits completely in $F$. Then $\mathfrak{X}^+(E/F_\cyc)$ has no nontrivial finite $\Lambda(\Gamma)$-submodules.
\end{enumerate}
\end{theorem}
\begin{remark}
More generally, suppose $E/F$ is an elliptic curve satisfying the following conditions (cf. Theorem 4.8 of \cite{LeiSujatha}):
\begin{enumerate}
    \item $E$ has supersingular reduction at all primes of $F$ above $p$.
    \item The map $\lambda_\cyc^\pm$ from (\ref{eq:lambdaplusminus}) is surjective.
    \item $H^2(F_S/F,E_{p^\infty})=0$.
    \item $H^1(\Gamma, \Sel_p^\pm(E/F_\cyc))=0$.
\end{enumerate}
Then, $\mathfrak{X}^\pm(E/F_\cyc)$ has no nontrivial finite $\Lambda(\Gamma)$-submodules.
\end{remark}
\begin{remark}
See Theorem \ref{Billotnopseudonull} for a similar result of Billot on elliptic curves with complex multiplication.
\end{remark}

\vs

Recall the Definition \ref{def:Eulerchar} of the Euler characteristic. Let $\Gamma=\Gal(F_\cyc/F)$ and thus, $\Gamma$ is topologically isomorphic to $\bZ_p$. Hence $\Gamma$ has $p$-cohomological dimension equal to $1$. Therefore, if $\chi(\Gamma, \Sel_p^\pm(E/F_\cyc))$ is defined, i.e. $H^i(\Gamma, \Sel_p^\pm(E/F_\cyc)$ is finite for $i=0$ and $1$, then 
\[
\chi(\Gamma, \Sel_p^\pm(E/F_\cyc))=\frac{\# H^0(\Gamma, \Sel_p^\pm(E/F_\cyc))}{\# H^1(\Gamma, \Sel_p^\pm(E/F_\cyc))}
\]

\vs

The following Theorem is due to Kim (Theorem 1.2, \cite{kimdu2013plus}).
\begin{theorem}\label{thm:cycEulerchar}
Let $E/F$ be an elliptic curve over a number field $F$ and let $\Gamma=\Gal(F_\cyc/F)$. Suppose $E$ has supersingular reduction at every prime above $p$. Furthermore, suppose that $\Sel_p(E/F)$ is finite. Then, $\chi(\Gamma, \Sel_p^\pm(E/F_\cyc))$ is well-defined, and up to a p-adic unit, is equal to
$$
\chi(\Gamma, \Sel_p^\pm(E/F_\cyc))= \#\Sel_p(E/F) \times \prod_\nu c_\nu
$$
where $\nu$ runs over all primes of $F$, $c_\nu$ is the Tamagawa number of $E$ at prime $\nu$ and $c_\nu=[E(F_\nu):E_0(F_\nu)]$. Here, $E_0(F_\nu)$ denotes the subgroup of points of $E(F_\nu)$ with non-singular reduction.
\end{theorem}
Finally, we mention a result of  analogous to the above result. 
\begin{remark}
Ahmed and Lim in \cite{ahmed_lim} extend the above result to compute the Euler character of the signed Selmer groups of elliptic curves with \textbf{mixed reduction type} (not necessarily supersingular reduction at all primes above $p$) over the cyclotomic $\bZ_p$-extension (cf. \cite[Theorem 1.1]{ahmed_lim}).
\end{remark}
\subsection{$\bZ_p^2$-extension case}
We work under the settings of Section \ref{sec:sec1.3} and follow the results of Lei and Sujatha \cite{LeiSujatha}. Let $H=\Gal(F_\infty/F_\cyc)\cong\bZ_p$ and $G=\Gal(F_\infty/F)$. In particular, note that:
\[
G/H=\Gal(F_\cyc/F)=\Gamma\cong\bZ_p.
\]
Recall the left exact sequence (\ref{eq:lambdaplusminus}). Let us consider the following commutative diagram (cf. diagram 5.2 of \cite{LeiSujatha}):
\begin{equation}\label{eq:diag}
\begin{tikzcd}
 0 \arrow[r] & \Sel_p^\pm(E/F_\infty)^H \arrow[r] & H^1(F_S/F_\infty,E_{p^\infty})^H \arrow[r, "\lambda_\infty^{H,\pm}"] &\displaystyle\oplus_{\nu \in S} J_\nu^\pm(E/F_\infty)^H\\
 0 \arrow[r] & \Sel_p^\pm(E/F_\cyc) \arrow[r] \arrow[u, "\alpha_\infty^\pm"] & H^1(F_S/F_\cyc,E_{p^\infty}) \arrow[r, "\lambda_\cyc^{\pm}"] \arrow[u, "\beta_\infty^\pm"] & \displaystyle\oplus_{\nu \in S} J_\nu^\pm(E/F_\cyc) \arrow[u, "\gamma_\infty^\pm=\oplus\gamma_{\nu,\infty}^\pm"].
 \end{tikzcd}
\end{equation}
Then, Lei and Sujatha proved the following results (cf. Section 5.3 of \cite{LeiSujatha}).
\begin{theorem}
The following statements hold true for the maps appearing in (\ref{eq:diag}):
\begin{enumerate}
    \item If $\lambda_{\cyc}^\pm$ is surjective, then $\alpha_\infty^\pm$ is \textit{injective} and $\T{coker}(\alpha_\infty^\pm) \cong \ker(\gamma_\infty^\pm)$. 
    \item The map $\beta_\infty^\pm$ is an isomorphism.
    \item The map $\gamma_{\nu,\infty}^\pm$ is an isomorphism for all $\nu \in \Sss_{p,F}$ and for all $\nu \in S$ with $\nu\not| p$.
    \item The map $\gamma_{\nu,\infty}^\pm$ is surjective and  $\ker(\gamma_{\nu,\infty}^\pm)$ is $\Lambda(\Gamma)$-cotorsion\footnote{A $\Lambda(G)$-module $M$ is a cotorsion if its Pontryagin dual is a torsion $\Lambda(G)$-module.} for all $\nu \in S \backslash \Sss_{p,F}$.
\end{enumerate}
\end{theorem}
\begin{remark}
In particular, note that if $E$ has good supersingular redaction at all primes above $p$, then all the vertical maps in (\ref{eq:diag}) become isomorphism maps. 
\end{remark}
Recall that Theorem \ref{thm:surjective} implies that if $\mathfrak{X}^{\pm}(E/F_\cyc)$ is torsion, then $\lambda_\cyc^\pm$ is surjective. Furthermore, Lei and Sujatha proved (cf. \cite{LeiSujatha}, Proposition 5.11):
\begin{theorem}
Suppose $\mathfrak{X}^{\pm}(E/F_\cyc)$ is torsion. Then the map $\lambda_\infty^{H,\pm}$ in (\ref{eq:diag}) is surjective.
\end{theorem}
\begin{remark}
We have a similar statement as Theorem \ref{thm:surjective} for the map
\begin{equation*}
    \lambda_\infty^\pm \ : \ H^1(F_S/F_\infty, E_{p^\infty}) \to \displaystyle\oplus_{\nu \in S} J_\nu^\pm(E/F_\infty).
\end{equation*}
Namely, $\mathfrak{X}^{\pm}(E/F_\infty)$ is $\Lambda(G)$-torsion exactly when $\lambda_\infty^\pm$ is surjective and $H^2(F_S/F_\infty, E_{p^\infty})=0$ (cf. proof of Corollary 5.12 in \cite{LeiSujatha}).
\end{remark}
\begin{remark}
Using Kim's control Theorem (cf. Theorem 1.1 of \cite{Kim}), Lei and Sujatha proved that (cf. \cite{LeiSujatha}, Corollary 5.12 and Proportion 5.14):
\begin{enumerate}[i)]
    \item If we assume $\Sel_p(E/F)$ is finite, then $H^1(H,\Sel_p^\pm(E/F_\infty))=0$.
    \item If we assume $\Sel_p(E/F)$ is finite and $E$ has supersingular reduction at every prime above $p$, then $H^i(G,\Sel_p^\pm(E/F_\infty))=0$ for $i \geq 1$.
\end{enumerate}
\end{remark}
The following Theorem \cite[Theorem 5.15]{LeiSujatha} extends Theorem \ref{thm:cycEulerchar} to $G$-Euler characteristic of  $\Sel_p^\pm(E/F_\infty)$.
\begin{theorem}\label{Thm:LeiSujEuler}
Let $E/F$ be an elliptic curve with  supersingular reduction at every prime above $p$. Furthermore, suppose $\Sel_p(E/F)$ is finite. Then, $\chi(G, \Sel_p^\pm(E/F_\infty))$ is well-defined, and 
$$
\chi(G, \Sel_p^\pm(E/F_\infty))=\chi(\Gamma, \Sel_p^\pm(E/F_\cyc)).
$$
\end{theorem}
\section{Conjecture A and $\mu$-invariant of torsion part of dual Selmer group}\label{Sec:ConjeAimplies}
In this section, we would like to show that Conjecture A implies that the torsion part of the $p$-Selmer group have $\mu$-invariant zero. To see this, we need to relate the fine Selmer group and the torsion part of the Selmer group.

\subsection{Results of Billot and Wingberg}
Suppose the following statements:
\begin{enumerate}[i)]
\item $E/F$ is an elliptic curve with complex multiplication by an imaginary quadratic field $K$ and $K \subset F$.
\item $  E_p \subset E(F)$.
\item $p$ is an odd prime number inert in $K$. 
\item $E$ has good reduction on all places of $F$ over $p$.
\end{enumerate}
Let $F_\infty=F(E_{p^\infty})$. Then, $\Gal(F_\infty/F)$ is topologically isomorphic to $\bZ_p^2$ \cite[Section 3]{Billot}. Let $G=\Gal(F_\infty/F)$ and let $\Lambda(G)$ be the Iwasawa algebra of $G$. Furthermore, we assume:
\begin{enumerate}[v)]
\item The Leopoldt's Conjecture\footnote{See Section 5.5 and 13.1 of \cite{washington1997introduction} for more information on Leopoldt's Conjecture.} holds for all the intermediate fields in the $\bZ_p^2$-extension $F_\infty$ of $F$.
\end{enumerate}
In the case where $E$ has good supersingular reduction for some of the primes above $p$, we expect the Pontryagin dual of the Selmer group $\mathfrak{X}(E/F_\infty)$ to have a positive rank. On the other hand, assuming conditions (i) to (iv), Billot proved that the dual of $\Sel^0(E/F_\infty)$ is a torsion $\Lambda(G)$-module (cf. Proposition 3.8 of \cite{Billot}). Let us denote the $\Lambda(G)$-torsion submodule of $\mathfrak{X}(E/F_\infty)$ by $T_{\Lambda(G)} (\mathfrak{X}(E/F_\infty))$. We have the following results by Billot.
\begin{theorem}\label{Thm:Billot}[\cite{Billot}, Theorem 3.23]
Given the above assumptions (i) to (v), the $\Lambda(G)$-torsion submodule of Pontryagin dual of the Selmer group $T_{\Lambda(G)} (\mathfrak{X}(E/F_\infty))$ is pseudo-isomorphic to the Pontryagin dual of the fine Selmer group $\mathfrak{X}^0(E/F_\infty)$.
\end{theorem}
\begin{theorem}[\cite{Billot}, Proposition 3.26]\label{Billotnopseudonull}
Suppose the assumptions (i) to (v) hold true.
\begin{enumerate}
    \item $\mathfrak{X}(E/F_\infty)$ has no non-trivial pseudo-null $\Lambda(G)$-submodule.
    \item Let $K_\infty$ be a $\bZ_p$-extension of $F$ contained in $F_\infty$ and let $H=\Gal(K_\infty/F)$. Suppose 
    \[
    \T{rank}_{\Lambda(G)}(\mathfrak{X}(E/F_\infty))=\T{rank}_{\Lambda(H)}(\mathfrak{X}(E/K_\infty)).
    \]
    Then, $\mathfrak{X}(E/K_\infty)$ has no non-trivial pseudo-null $\Lambda(H)$-submodule.
\end{enumerate}
\end{theorem}
 Now, suppose $F$ is a number field and $K_\infty/F$ is a $\bZ_p$-extension of $F$. The Tate-Shafarevich group $\Sha(E/F)$ of elliptic curve $E/F$ is defined by the exact sequence 
$$
0 \to \Sha(E/F) \to H^1(F,E) \to \displaystyle\bigoplus_\omega H^1(F_\omega,E),
$$
where $\omega$ runs over all archimedean and non-archimedean places of $F$. We let $\Sha(E/F)_{p^i}$ to be the $p^i$-th torsion subgroup of $\Sha(E/F)$ and $\Sha(E/F)_{p^\infty}=\varinjlim_i \Sha(E/F)_{p^i}$. Let $F_n/F$ be the $n$-th layer of the cyclotomic $\bZ_p$-extension $F_\cyc/F$. Here, we denote $\Sha(E/F_n)_{p^\infty}$ by $\Sha_n$. The following Theorem by Wingberg (\cite{Wingberg}, Corollary 2.5) is a generalization of Billot's result.
\begin{theorem}\label{Thm:Wingberg}
Suppose $E/F$ is an elliptic curve with supersingular reduction for every prime of $F$ over $p$. Assume that $\Sha_n$ is finite for all $n$ and the Pontryagin dual of the $E(F_\cyc) \otimes \bQ_p/\bZ_p$ is $\Lambda(\Gamma)$-torsion. Then, 
$$
T_{\Lambda(\Gamma)} (\mathfrak{X}(E/F_\cyc)) \sim \mathfrak{X}^0(E/F_\cyc).
$$\end{theorem}
\begin{remark}
Assuming Conjecture A, we get that $T_{\Lambda(\Gamma)} (\mathfrak{X}(E/F_\cyc))$ has $\mu$-invariant equal to zero.
\end{remark}
\begin{remark}
Wingberg shows that the dual fine Selmer group can be identified with the adjoint of $M$, where $M$ is the torsion submodule of the dual Selmer group. The adjoint of $M$ is defined as $\Ext^1_{\Lambda(G)}(M,\Lambda(G))$ (cf. page 473 of \cite{Wingberg}). This adjoint is pseudo-isomorphic to $M$, as $M$ is a torsion $\Lambda(G)$-module (cf. \cite{perrin1984arithmetique}, Chapter 1, Section 2.2).
\end{remark}
\begin{remark}
In \cite{Matar}, Ahmed Matar gives a Galois theoretic proof of Wingberg's result (Theorem \ref{Thm:Wingberg}) with a slightly different hypothesis. See Theorem 1.1 of \cite{Matar}.
\end{remark}
\subsection{Numerical examples}\label{Numerical examples}
In this section, we give some evidence for Conjecture A. We use the results of C. Wuthrich \cite{Wuthrich}.

\begin{theorem}\label{thm:isogeny}[\cite{Wuthrich}, Proposition 8.1]
Suppose $p$ is an odd prime and $E/\bQ$ is an elliptic curve which admits an isogeny over $\bQ$ of degree $p$. Then, the $\mu$-invariant of the Pontryagin dual of the fine Selmer group $\mathfrak{X}^0(E/\bQ_\cyc)$ is zero as a $\Lambda(\Gamma)$-module.
\end{theorem}
Recall that, by Theorem \ref{Kobayashi 7.2},  $\mathfrak{X}^0(E/\bQ_\cyc)$ is a finitely generated $\Lambda(\Gamma)$-module. Therefore, if $E/\bQ$ satisfies the hypothesis of Theorem \ref{thm:isogeny}, then Conjecture A is satisfied for $E/\bQ$.

\vs

Here we summarize several numerical examples of elliptic curves $E/\bQ$ with a description of their fine Selmer groups for specific primes $p$ discussed in \cite{Wuthrich}, in Sections 9, 10, and 11.\\

\begin{center}
\begin{tabular}{ |c|c|c|}
\hline
\begin{tabular}{@{}c@{}} Cremona label \\ of $E/\bQ$ \end{tabular}  & \begin{tabular}{@{}c@{}}Mordell--Weil\\rank of $E$\end{tabular} &\begin{tabular}{@{}c@{}} Fine Selmer groups \\  $\Sel_p^0(E/\bQ_\cyc)$\end{tabular} \\
\hline
11A1    & 0 &  \begin{tabular}{@{}c@{}} Trivial for all odd primes $p\neq5$. \\ Finite (but none-trivial) for $p=5$ \end{tabular} \\
\hline
11A2    & 0 &  Trivial for all odd primes\\
\hline
11A3    & 0 &  Trivial for all odd primes \\
\hline
182D1   & 0 &  Conjecture A holds for $p=5$\\
\hline
37A1   & 1  &  Finite for all odd primes $p<1000$\\
\hline
53A1    & 1 &  Finite for $p=3$\\
\hline
5692A1  & 2 &  Conjecture A holds for $p=3$\\
\hline
\end{tabular}
\end{center}

\vs

\begin{remark}
Note that the Pontryagin dual of a finite (resp. trivial) group is again finite (resp. trivial). Therefore, when the fine Selmer $\Sel_p^0(E/\bQ_\cyc)$ is finite, its Pontryagin dual $\mathfrak{X}^0(E/\bQ_\cyc)$ is finite (and so pseudo-null). In particular, Conjecture A holds here.
\end{remark}
\begin{remark}
The first three elliptic curves (11A1, 11A2, and 11A3) in the above table represent all the elliptic curves $E/\bQ$ of conductor $11$ upto isogeny.
\end{remark}

\section{Signed Selmer groups and $p$-adic Hodge Theory}\label{Sec:pAdicHodge}

In this section, we revisit signed Selmer groups using Fontaine's $p$-adic Hodge Theory. These definitions of signed Selmer groups are due to Lei, Loeffler and Zerbes (see \cite{LeiLoefflerZerbes}, \cite{LeiZerbes}, \cite{Lei_IwasawaTheory}).  They defined the "local conditions" (cf. \eqref{eq:localcondtions}) on these Selmer groups  using Fontaine's rings  $\mathbb{E}$ and $\widetilde{\mathbb{E}}$ (cf. Section \ref{sec:perfecthull}) over the extension $L(\mu_{p^\infty})$ where $L$ is a finite extension of $\bQ_p$. Recently, Scholze realized that the Fontaine's ring $\widetilde{\mathbb{E}}$ is the tilt of the perfectoid field $\widehat{L(\mu_{p^\infty})}$ (cf. Theorem \ref{Thm:Scholze}). This tilting construction works in general for any perfectoid field and Scholze also generalized Fontaine-Wintenberger Theorem \cite{FontaineWintenberger} for Galois groups of arbitrary perfectoid fields and their tilts. Therefore, a natural question is to ask whether the local conditions on the signed Selmer groups can be generalized using the tilt of any general perfectoid field and thereby extend the definition of signed Selmer groups for any general $p$-adic Lie extension (other than the cyclotomic extension, for which they are already defined by Kobayashi). We record this as an open question (cf. see Question \ref{ques:perfectoid}). This is the motivation behind this section and therefore we start by defining the Fontaine rings using the language of perfectoid spaces  of Scholze \cite{scholzeberkeley}.

\subsection{Perfectoid and Tilts}
Let $L$ be a finite extension of $\bQ_p$ with ring of integers $\cO_L$. Here, we will write $\Gamma_L$ for the Galois group $\Gal(L(\mu_{p^\infty})/L)$ which we identify with $\cO_L^*$ via the cyclotomic character $\chi_L$. Fix $\pi$ a uniformizer of $L$ normalized by $|\pi|=q^{-1}$ where $q=|k_L|$ , the cardinality of the residue field $k_L$ of $L$. Here, we will follow the exposition of \cite{schneider2017galois}. 
\begin{definition}
Any intermediate field $L \subset K \subset \bC_p$ is called perfectoid, if it satisfies the following conditions:
\begin{enumerate}
\item $K$ is complete,
\item The value group $|K^*|$ is dense in $\bR^*_{>0}$ , and 
\item $(\cO_K/p\cO_K)^p=\cO_K/p\cO_K$.
\end{enumerate}
\end{definition}
One can show that every element of the value group $|K^*|$ is a $p$th power.\\
From the field $K$ of characteristic $0$, we will construct a new field $K^\flat$ (the tilt of $K$) of characteristic $p$. We first choose an element $\varpi\in \mathfrak{m}_K$ such that $|\varpi| \geq |\pi|$ and define
\begin{align*}
\cO_{K^\flat}&:=
\varprojlim\left(
\cdots \xrightarrow{(\cdot)^q}  \cO_K/\varpi\cO_K  \xrightarrow{(\cdot)^q}  \cO_K/\varpi\cO_K \xrightarrow{(\cdot)^q} \cdots \xrightarrow{(\cdot)^q} \cO_K/\varpi\cO_K \right)
\\
&=\{(\cdots,\alpha_i,\cdots,\alpha_1,\alpha_0) \in (\cO_K/\varpi\cO_K)^{\bN_0} \ : \ \alpha_{i+1}^q=\alpha_i \ \T{ for all } i \geq 0\}.
\end{align*}
Let $K^\flat$ be the fraction field of the integral domain $\cO_{K^\flat}$.
Let $k$ be the residue field of $K$.
The $k$-algebra $\cO_{K^\flat}$ is perfect and one can further show that, for any element $\alpha=(\cdots, \alpha_i,\cdots, \alpha_0)\in \cO_{K^\flat}$, if we choose $a_i \in \cO_{K}$ such that $a_i$ mod $\varpi \cO_K=\alpha_i$, then the limit
\[
\alpha^\sharp:=\varinjlim_{i \to \infty} a_i^{q^i} \in \cO_K
\]
exists and is independent of the choice of lifts $a_i$'s. This gives the (untilt) map
\begin{align*}
    \cO_{K^\flat} &\to \cO_K \\
    \alpha &\to \alpha^{\sharp}
\end{align*}
which is a well-defined multiplicative map such that $\alpha^\sharp$ mod $\varpi\cO_K=\alpha_0$. We have a multiplicative bijection 
\[
\cO_{K^\flat} \xrightarrow{\sim} \varprojlim_{(\cdot)^q} \cO_K
\]
given by 
\[
\alpha \mapsto (\cdots, (\alpha^{\frac{1}{q^i}})^{\sharp},\cdots,\alpha^\sharp)
\]
which shows that $\cO_{K^\flat}$ is independent of the choice of the element $\varpi$.
\subsection{Values on Tilt}
The map 
\begin{align*}
| \cdot |_\flat : \cO_{K^\flat} &\to \bR_{\geq 0}\\
 \alpha &\to |\alpha^\sharp|
\end{align*}
is a non-archimedean absolute value such that
\begin{enumerate}[i)]
\item $|\cO_{K^\flat}|_\flat=|\cO_K|$; 
\item $\mathfrak{m}_{K^\flat}:=\{\alpha \in \cO_{K^\flat} \ : |\alpha|_{K^\flat} <1 \} $ is the unique maximal ideal in $\cO_{K^\flat}$ and we have $\cO_{K^\flat}/\mathfrak{m}_{K^\flat}\cong\cO_{K}/\mathfrak{m}_{K}$;
\item Let $\varpi^\flat \in \cO_{K^\flat}$ be any element such that $|\varpi^\flat|_{\flat}=|\varpi|$, then we have an isomorphism $\cO_{K^\flat}/\varpi^\flat\cO_{K^\flat} \cong \cO_K/\varpi\cO_K$.
\end{enumerate}
 Since $K^\flat$ is the fraction field of the integral domain $\cO_{K^\flat}$, the norm $| \cdot |_\flat$ extends to a non-archimedean absolute value on $K^\flat$. One can show that the value group of $K^\flat$ under the norm $| \cdot |_\flat$ is same as the value group of $K$ under $| \cdot |$. Under the norm $| \cdot |_\flat$,  $\cO_{K^\flat}$ is the ring of integers of $K^\flat$. Further, $K^\flat$ with $| \cdot |_\flat $ is a perfect and complete non-archimedean field of characteristic $p$. When $K=\bC_p$, Lemma 1.4.10 of \cite{schneider2017galois} shows that its tilt $\bC_p^\flat$ is algebraically closed.

\vs

An important example of perfectoid fields is the completion $\widehat{L(\mu_{p^\infty})}$ of the field $L(\mu_{p^\infty})$ (cf. Proposition 1.4.12 of \cite{schneider2017galois}). We have a natural continuous action (cf. Lemma 1.4.13 of \cite{schneider2017galois}) of $\Gal(\overline{\bQ}_p/L)$ on $\bC_p^\flat$ which is continuous with respect to the topology induced by $| \cdot |_\flat$.

\vs

The subgroup $\Gal(\overline\bQ_p/L(\mu_{p^\infty}))$ fixes $\widehat{L(\mu_{p^\infty})}^\flat \subseteq \bC_p^\flat$ and hence we have a continuous action of $\Gamma_L=Gal(L(\mu_{p^\infty})/L)$ on the tilt $\widehat{L(\mu_{p^\infty})}^\flat$. Let $\varepsilon=(\cdots, \varepsilon_i, \cdots, \varepsilon_1, \varepsilon_0) \in \widehat{L(\mu_{p^\infty})}^\flat$ such that $\varepsilon_0=1$, $\varepsilon_1\neq 1$ and let $X=\varepsilon-1$. Then, we have a well defined embedding of fields (cf. p. 50 of \cite{schneider2017galois}) 
\begin{equation}\label{eq:E_L}
k((X)) \xrightarrow{\theta} \widehat{L(\mu_{p^\infty})}^\flat,
\end{equation}
where $k$ is the common residue field of $L(\mu_{p^\infty}), \widehat{L(\mu_{p^\infty})}$, and $ \widehat{L(\mu_{p^\infty})}^\flat$ (cf. Proposition 1.3.12(i) and Lemma 1.4.6(iv) of \cite{schneider2017galois}). We denote the image of $\theta$ by $\mathbb{E}_L$. Then, $(\mathbb{E}_L,| \cdot |_\flat)$ is a complete non-archimedean discretely valued field with residue field $k$ which is preserved by the $\Gamma_L$-action on $\widehat{L(\mu_{p^\infty})}^\flat$.\\
\subsection{Perfect hull}\label{sec:perfecthull}
Suppose $E$ is any field of characteristic $p>0$ and let $\overline{E}/E$ be an algebraic closure. Then, we define 
\[
E^{perf}:=\{ a\in \overline{E} \ : \ a^{p^m} \in E \text{  for some $m\geq 0$}  \},
\]
which is called the \textbf{perfection (or perfect hull)} of $E$. In particular, $E^{perf}$ is the smallest intermediate field of $\overline{E}/E$ which is perfect.
Here are two most crucial Theorems whose proofs can be found in \cite{schneider2017galois} (Proposition 1.4.17 and Proposition 1.4.27).
\begin{theorem}[Scholze]\label{Thm:Scholze}
Let $\mathbb{E}_L$ be the image of the injective map $\theta$ in \eqref{eq:E_L}. Then, we have:
\begin{enumerate}
\item $\widehat{\mathbb{E}_L^{perf}}=\widehat{L(\mu_{p^\infty})}^\flat$,
\item $ \widehat{\mathbb{E}_L^{sep}}=\bC_p^\flat$.
\end{enumerate}
\end{theorem}

\vs

For simplicity, we let $\mathbb{E}:=\mathbb{E}_L^{sep}$ and $\widetilde{\mathbb{E}}:=\widehat{\mathbb{E}_L^{perf}}$. We can identify $\mathcal{H}_L=\Gal(\overline{L}/L(\mu_{p^\infty}))$ with $\Gal(\mathbb{E}/\mathbb{E}_L)$; this is Fontaine-Wintenberger's Theorem \cite{FontaineWintenberger}. (This now holds for more general perfectoid fields by the work of Scholze).  By continuity, the action of $\Gal(\overline{L}/L)$ on $\mathbb{E}$ extends uniquely to an action on $\widetilde{\mathbb{E}}$ and we let $\widetilde{\bE}_L:=\widetilde{\bE}^{\cH_L}$. Let $\widetilde{\bA}$ be the ring of Witt vectors of $\widetilde{\bE}$, $\widetilde{\bB}=\widetilde{\bA}[p^{-1}]$.

\vs

Recall from \eqref{eq:E_L}, that $\bE_L\cong k((\varepsilon-1))$ and let $\bE^+_L\cong k[[\varepsilon-1]]$ be its valuation ring. Let $\bA_L^+$ be a complete regular local ring of dimension $2$ such that $\bA_L^+/(p) \cong \bE_L^+$. Such a lift exists and it is unique (\cite{scholl2006higher} p. 699 and \cite{berger2004introduction}, Section III.1.1). Let $\bA_L$ be the $p$-adic completion of $\bA_L^+[\frac{1}{\varepsilon-1}]$. Define $\bB_L^+:=\bA_L^+[\frac{1}{p}]$ and $\bB_L:=\bA_L[\frac{1}{p}]$.
	
\vs
	
The rings $\bA_L$, $\bA_L^+$, $\bB_L$, and $\bB_L^+$ are endowed with an action of $\Gamma$, Frobenius $\varphi$, and its left inverse $\psi$. When $L$ is unramified, we have simple formulae for $\Gamma, \ \varphi$, and $\psi$-action. Write $X:=\varepsilon-1$ and $\gamma \in \Gamma$, then we have (cf. \cite[ 2.3]{loefflerZerbes}),
\begin{align*}
		\gamma(X)&= (1+X)^{\chi_L(\gamma)}-1, \\
		\varphi(X)&=(1+X)^p-1,\\
		\varphi \circ \psi (f(X))&=\frac{1}{p}\ds\sum_{\zeta^p=1}f\left(\zeta(1+X)-1\right).
\end{align*}
Let $\bB$ be the $p$-adic completion of the maximal unramified extension of $\bB_{\bQ_p}=\bA_{\bQ_p}[p^{-1}]$ in $\widetilde{\bB}=W(\widetilde{\bE})[\frac{1}{p}]=\widetilde{\bA}[\frac{1}{p}]$ and let $\bA=\bB \cap \widetilde{\bA}$.
	
\vs
	
These rings are also stable under Frobenius and the absolute Galois group $\cG_{\bQ_p}$ of $\bQ_p$ (under this identification, $\bA_L=\bA^{\cH_L}$). We write $\T{Rep}_{\bZ_{p}}(\cG_{\bQ_p})$ (resp. $\T{Rep}_{\bQ_{p}}(\cG_{\bQ_p})$) for the category of finitely generated $\bZ_p$-modules (resp. finite dimensional $\bQ_p$-vector space) with a canonical action of $\cG_{\bQ_p}$.
	
\vs
	
For a $p$-adic representation $T \in \T{Rep}_{\bZ_p}(\cG_{\bQ_p})$ (resp. $V\in \T{Rep}_{\bQ_p}(\cG_{\bQ_p})$), we define a free finitely generated module over $\bA_L$ of rank equal to  $\T{rank}_{\bZ_p}(T)$ (resp. a finite dimensional vector space over $\bB_L$ of dimension equal to $\dim_{\bQ_{p}}(V)$) by $D_L(T)=(\bA \ot_{\bZ_p}T)^{\cH_L}$ (resp. $D_L(V)=(\bB \ot_{\bQ_p}V)^{\cH_L}$). These are equipped with a commuting semi-linear action of $\varphi$ and $\Gamma_L$. Further, we can define the left inverse $\psi$ of $\varphi$ on $\bA$ which then extends to a left inverse $\psi$ of $\varphi$ on $D_L(V)$ and $D_L(T)$ (cf. \cite{LeiZerbes}).

\subsection{Iwasawa Cohomology}

$T \in \T{Rep}_{\bZ_p}(\cG_L)$ and let $L_\infty$ be a $p$-adic Lie extension of $L$. We define
\[
H^i_{\T{Iw}}(L_\infty,T):=\ds\varprojlim_K H^i(K,T),
\]
where $K$ varies over the finite extension of $L$ contained in $L_\infty$ and the inverse limit is taken with respect to the corestriction maps.

\vs

If $V=\bQ_p \ot_{\bZ_p} T$, then we write
\[
H^i_{\T{Iw}}(L_\infty,V)=\bQ_p \ot_{\bZ_p} H^i_{\T{Iw}}(L_\infty,T),
\]
which is independent of the choice of the lattice $T \subset V$. These are finitely generated $\Lambda(G)$-modules, where $G=\Gal(L_\infty/L)$ (cf. Theorem A.2 of \cite{loefflerZerbes}). 

\subsection{Fontaine Isomorphism}
In the case when $L_\infty=L(\mu_\infty)$, we have a canonical isomorphism of $\Lambda(\Gamma_L)$-modules:
\[
h^1_{\T{Iw}, T}: \bD_L(T)^{\psi=1} \xrightarrow{\cong} H^1_{\T{Iw}}(L(\mu_{p^\infty}),T).
\]
Further, if $L$ is an unramified extension of $\bQ_p$ and $V=\bQ_p \ot_{\bZ_p} T$ is a crystalline representation of $\Gal(\overline{L}/L)$ with non-negative Hodge-Tate weights and no quotient isomorphic to $\bQ_p$, then we have
\[
h^1_{\T{Iw},T}: \bN_L(T)^{\psi=1} \xrightarrow{\cong} H^1_{\T{Iw}}\left(L(\mu_{p^\infty}),T\right)
\]
where $\bN_L(T)$ is the Wach module associated to $T$, due to Berger (the reader should see Section 2.7 of \cite{loefflerZerbes} as a reference for these results).
\subsection{Iwasawa Cohomology and $\mu$-invariant}
Recall the structure Theorem of finitely generated modules over the Iwasawa algebra (see Theorem \ref{thm:structurethm} for more details). By Shapiro's Lemma, we know that for $L_\infty=L(\mu_{p^\infty})$, we have
\[
H^i_{\T{Iw}}(L_\infty,V)=H^i\left(\cG_L,\cO_L[[\Gamma_L]] \ot_{\cO_L} V\right),
\]
where the right hand side is the usual Galois cohomology with continuous cochains and $\cG_L=\Gal(\overline{L}/L))$. We have the following fact from \cite[Lemma 5.12]{SchVenj}.
\begin{lemma}
\begin{enumerate}[i)]
    \item $H^i_{\T{Iw}}(L_\infty,V)=0$ for $i\not= 1,2$.
    \item $H^2_{\T{Iw}}(L_\infty,V)$ is a finitely generated $\cO_L$-module (and therefore has $\mu$-invariant zero).
\end{enumerate}
\end{lemma}
\subsection{$p$-adic Hodge Theory and Signed Selmer groups}
In this Section, we follow the exposition in \cite{LeiZerbes} and discuss how we can define signed Selmer groups using $p$-adic Hodge theory.
\subsubsection{Good supersingular elliptic curves} Suppose $E$ is an elliptic curve over $\bQ$. Let $p\geq 3$ and $a_p=0$. Let $T:=T_p(E)$ be the $p$-adic Tate module of $E$ and write $V=T_p(E) \ot_{\bZ_p} \bQ_p$ which is a crystalline representation of $\cG_{\bQ_p}$ with Hodge-Tate weights $0, \ 1$. Let $\{v_{+},v_-\}$ be a basis of $\bD_{cris}(V)=(\bB_{cris} \ot_{\bQ_p} V)^{\cG_{\bQ_p}}$ (where  $\bB_{cris}$ is one of the Fontaine's period rings) such that the matrix of $\varphi$ on $\bD_{cris}(V)$ in this basis is
\[
\begin{bmatrix}
0 & -1\\
p & 0
\end{bmatrix}.
\]
Let $q=\varphi(\pi)/\pi$ where $\pi:=[\varepsilon]-1$. Define 
\[
\log_+(1+\pi)=\ds\prod_{i\geq 0} \frac{\varphi^{2i}(q)}{p},
\]
and
\[
\log_-(1+\pi)=\ds\prod_{i\geq 0} \frac{\varphi^{2i+1}(q)}{p},
\]
	where $\varepsilon=(\varepsilon^{(i)})$ was the fixed element in $\widetilde{\bE}$ such that $\varepsilon^{(0)}=1$, $\varepsilon^{(1)}\not=1$, and $[\varepsilon]$ is the Teichmüller lift of $\varepsilon \in \widetilde{\bE}$. Suppose 
	\[
	\begin{bmatrix}
	n_+\\
	n_-
	\end{bmatrix}
	=M
	\begin{bmatrix}
	v_+\\
	v_-
	\end{bmatrix},
	\]
where
\[
M=	
\begin{bmatrix}
\log_+(1+\pi) & 0\\
0 & \log_-(1+\pi)
\end{bmatrix}.
\]
Then, $n_+$ and $n_-$ from a basis of the Wach module $\bN(T)$ (cf. Section 3.1 of \cite{LeiZerbes}, see also \cite{LeiLoefflerZerbes_WachModules}).

\vs

Suppose $L$ is a finite extension of $\bQ_p$. Then, any element 
\[
x \in \bN_L(T)^{\psi=1}\cong \bD_{L}(T)^{\psi=1}\cong H^1_{Iw}(L(\mu_{p^\infty}),T)
\]
can be written as 
\[x=x_+v_++x_-v_-=x_+^{\prime}n_++x_-^{\prime}n_-,\]
where $x_+=x_+^\prime \log_+(1+\pi)$ and $x_-=x_-^\prime \log_-(1+\pi)$. Then, for $i=\pm$, we define 
\[
H^1_{Iw}(L(\mu_{p^\infty}),T)^i:=\{x\in\bN_L(T)^{\psi=1} \ : \ \varphi(x_i)= - p \psi(x_i)  \}.
\]
For $n\geq 1$, we define $H^1(L_n,T)^i$ to be the image of $H^1_{Iw}(L(\mu_{p^\infty}),T)^i$ under the natural projection map
\[
H^1_{Iw}(L(\mu_{p^\infty}),T) \to H^1(L_n,T),
\]
where $L_n= L(\mu_{p^n})$.

\vs

In \cite{LeiLoefflerZerbes}, Section 5.2.1, the authors define two Coleman maps
\[
\T{Col}_i=H^1_{Iw}(L(\mu_{p^\infty}),T) \to \Lambda(\Gamma),
\]
and by remark 3.2 of \cite{LeiZerbes}, we know that
\[
H^1(L(\mu_{p^\infty}),T)^i=H^1_{Iw}(L(\mu_{p^\infty}),T) \cap \Ker(\T{Col}_i).
\]
Let $H^1_{f}(L_n,E_{p^{\infty}})^\pm$ be the exact annihilator of $H^1(L_n,T)^\pm$ under the Pontryagin duality
\[
[\sim,\sim]: H^1(L_n,T) \times H^1(L_n,E_{\infty}) \to \bQ_p/\bZ_p.
\]
 Now, suppose $F$ is a finite extension of $\bQ$ and the elliptic curve $E$ is defined over $\bQ$. Let $F_{\nu,n}:=F_{\nu}(\mu_{p^n})$ for a prime $\nu$ of $F$. Let $S$ be a finite set of primes of $\bQ$ containing $p$ and all the primes where $E$ has bad reduction and the infinite primes. For  $v \in S$, define the "new" local conditions
\begin{equation}\label{eq:localcondtions}
J_v^\pm(F_n) = \displaystyle\bigoplus_{w_n \mid v} \frac{H^1 (F_{w_{n},n},E_{p^\infty})}{H^1_{f}(F_{w_{n},n},E_{p^\infty})^\pm},
\end{equation}

for $w_n$ prime of $F_n$ over $\nu$.

We write $J_v^\pm(F_\cyc)=\varinjlim_nJ_v^\pm(F_n)$. We define
\begin{align*}
\Sel_p^{\pm}(E/F_n) &:= \ker \left(  H^1(F_S/F_n,E_{p^\infty}) \rightarrow \displaystyle\bigoplus_{v \in S} J_v^\pm(E/F_n)  \right);\\
\Sel_p^\pm(E/F_\cyc)&:=\varinjlim_n\Sel^\pm(E/F_n).
\end{align*}

Further, Lei shows that the "local conditions" $H^1_{f}(F_{w_{n},n},E_{p^\infty})^\pm$ in the definition of  $J_v^\pm(F_n)$ exactly coincides with $E^\pm(F_{w_{n},n}) \otimes \bQ_p/\bZ_p$ of Kobayashi (see the paragraph after the proof of Lemma $4.12$ of \cite{Lei_IwasawaTheory}). Therefore, this $p$-adic Hodge theory definition of Selmer groups is the usual Kobayashi's $\pm$ Selmer groups.

\section{Analogy with the ordinary reduction case}\label{Sec:Future}
Here, we list some of the main results proved in the case where we consider elliptic curves with good ordinary reduction over primes above $p$ and we can find analogous results for elliptic curves with supersingular reduction. The main take away here is that, many of the properties that the Selmer groups enjoy in the ordinary reduction case are observed for the signed Selmer groups when the reduction type is  supersingular.

\vs

Suppose $p\geq 5$ is an odd prime. We have the following Conjecture by Mazur (\cite{Mazur})\footnote{This Conjecture can be extended to other Galois extensions where the Galois group is a general $p$-adic Lie group.}.
\begin{conjecture}\label{conj:Mazure}
Suppose $E/F$ has good ordinary reduction at all primes above $p$. Then, $\mathfrak{X}(E/F_\cyc)$ is a finitely generated\footnote{It is not hard to show that, for any $\bZ_p$-extension $K_\infty/F$ of $F$, $\mathfrak{X}(E/K_\infty)$ is a finitely generated module over the Iwasawa algebra of $Gal(K_\infty/F)$ (cf. Lemma 2.4 of \cite{CoatesSujatha_book}).} $\Lambda(\Gamma)$-torsion module.
\end{conjecture}
\begin{remark}
This has been shown in several cases. For example, Mazur proved in \cite{Mazur} that, if $\Sel_p(E/F)$ is finite, then the above Conjecture holds. Furthermore, Kato in \cite{Kato} proved this Conjecture for the case where $E/\bQ$ is an elliptic curve defined over $\bQ$. 
\end{remark}
Often, the results proved for elliptic curves with good ordinary reduction motivate results in the supersingular reduction case like Theorem \ref{sspsuedonull}.
We have the following (cf. \cite{greenberg1997structure}).
\begin{theorem}
Suppose $E/F$ is an elliptic curve with good ordinary reduction at all primes above $p$. Assume $\mathfrak{X}(E/F_\cyc)$ is a finitely generated torsion  $\Lambda(\Gamma)$-module. Furthermore, assume $E(F)$ has no element of order $p$. Then, $\mathfrak{X}(E/F_\cyc)$ has no non-zero finite $\Lambda(\Gamma)$-submodule. 
\end{theorem}
The following Theorem (cf. \cite{perrin1992theorie}) is analogous to Theorem \ref{thm:surjective}.
\begin{theorem}
Suppose $p$ is an odd prime and $E/F$ is an elliptic curve with good ordinary reduction at all primes above $p$. Then, $\mathfrak{X}(E/F_\cyc)$ is a finitely generated torsion  $\Lambda(\Gamma)$-module if and only if the sequence
\begin{equation}
0 \to \Sel_p(E/F_\cyc) \to  H^1(F_S/F_\cyc,E_{p^\infty}) \xrightarrow{\lambda_\cyc} \ds\bigoplus_{\nu\in S} J_v(E/F_\cyc) 
\end{equation}
is short exact and $H^2(F_S/F_\cyc, E_{p^\infty})=0$. 
\end{theorem}

\vs

Let 
\[
\widetilde{F}_\infty:=F(E_{p^{\infty}})=\ds\cup_{n \geq 0}F(E_{p^{n}}).
\]
$ \widetilde{F}_\infty $ is often referred as the trivializing extension of $F$ (with respect to $E$ and prime $ p $). Let $\mathfrak{X}(E/\widetilde{F}_\infty)$ denote the Pontryagin dual of the Selmer group $\Sel_p(E/\widetilde{F}_\infty)$ and $ \widetilde{G}:=\Gal(\widetilde{F}_\infty/F) $. We mention the following result of Ochi-Venjakob (cf. \cite{ochi2002structure}, Theorem 5.1).
\begin{theorem}
Let $E/F$ be an elliptic curve without complex multiplication. Let $p \geq 5$ and suppose $ E $ has good ordinary reduction at all primes above $p$. Moreover, assume that $ \mathfrak{X}(E/\widetilde{F}_\infty) $ is a finitely generated $ \Lambda(\widetilde{G}) $-torsion module. Then, $ \mathfrak{X}(E/\widetilde{F}_\infty) $ has no non-zero pseudo-null $\Lambda(\widetilde{G})$-submodule.
\end{theorem}
\begin{remark}
The reader may also look at a similar Theorem by B. Perrin-Riou for the CM case (\cite{perrin1981groupe}, Theorem 2.4).
\end{remark}
Recall the definition \ref{def:Eulerchar} of $\Gamma$-Euler characteristic. The following Theorem is analogous to Theorem \ref{thm:cycEulerchar} (cf. Theorem 3.3 of \cite{CoatesSujatha_book}).
\begin{theorem}
Suppose $E/F$ has good ordinary reduction at all primes above $p$ and $\Sel_p(E/F)$ is finite. Then, $\chi(\Gamma, \Sel_p(E/F_\cyc))$ is well-defined is given by
$$
\chi(\Gamma, \Sel_p(E/F_\cyc))= \frac{\#\Sha(E/F)_{p^{\infty}}}{\#(E(F)_{p^{\infty}})^2} \times \ds\prod_\nu c^{(p)}_\nu
\times \ds\prod_{\nu\mid p}(d_\nu^{(p)})^{2}.
$$
Here, 
\begin{enumerate}
\item $\nu$ runs over all finite primes of $F$.
\item $c_\nu$ is the Tamagawa number of $E$ at prime $\nu$ and $c_\nu:=[E(F_\nu):E_0(F_\nu)]$, where $E_0(F_\nu)$ denotes the subgroup of points of $E(F_\nu)$ with non-singular reduction.
\item $ d_\nu:=\#(\tilde{E}_\nu(f_\nu)) $, where $f_\nu$ is the residue field of $F_\nu$  and $ \tilde{E}_\nu(f_\nu) $ is the reduction of $E$ modulo $\nu$ modulo $\nu$.
\item $c^{(p)}_\nu$ (resp. $ d_\nu^{(p)} $) is the largest power of $p$ dividing $ c_\nu$ (resp. $d_\nu  $).
\end{enumerate}
\end{theorem}

\vs

Now, suppose $E/F$ is an elliptic curve with complex multiplication by the ring of integers $\cO_F$. Suppose $F$ is an imaginary quadratic field with class number one. Let $ p $ be an odd prime such that it splits in $F$, $p\cO_F=\fp\overline{\fp}$ and suppose $E$ has good reduction at $\fp$ and $\overline{\fp}$. Since $E/F$ has complex multiplication, then $ \widetilde{G}=\Gal(\widetilde{F}_\infty/F) $ is topologically isomorphic to $\bZ_p^2$ and $ \Lambda(\widetilde{G}) \cong \bZ_p[[T_1,T_2]] $.  We have (cf. Proposition 3.2 of \cite{sujatha2011mu}):
\begin{theorem}
$\mathfrak{X}(E/\widetilde{F}_\infty)$ is a finitely generated torsion $ \Lambda(\widetilde{G}) $-module and $ \mu_{\widetilde{G}}(\mathfrak{X}(E/\widetilde{F}_\infty))=0 $.
\end{theorem}
To conclude, we pose two questions that remain open.
\begin{question}\label{ques:perfectoid}
Can we define the signed Selmer groups for general $p$-adic Lie extensions of a given number field?
\end{question}

\vs

In cases where signed Selmer groups are defined, many of the results that have been proved for the Selmer groups of elliptic curves with good ordinary reduction are expected to be true in the supersingular reduction for signed Selmer groups under suitable hypothesis. For example, we may ask:  
\begin{question}
Is it possible, under suitable conditions for an elliptic curve with supersingular reduction over primes above $p$, to prove that the $\mu$-invariants of the Pontryagin dual of the signed Selmer groups (also known as plus/minus $\mu$-invariants) vanish over the cyclotomic extension?
\end{question}

Answering these questions should involve some new ideas.

\section*{Acknowledgements}
We would like to thank Denis Benois and the organizers of Iwasawa 2019 conference for inviting the second author to speak at the conference. Thanks are also due to Sujatha Ramdorai for encouraging us to look at Conjecture A and applying it in the context of supersingular elliptic curves.

\bibliographystyle{alpha}
\bibliography{supersingular}
\end{document}